\documentclass[12pt,leqno]{article}
\usepackage{amsmath, amsthm, amsfonts, amssymb, color}
\usepackage{mathrsfs}

\theoremstyle{definition}

\newcommand{\scr}[1]{\mathscr #1}
\definecolor{wco}{rgb}{0.5,0.2,0.3}

\numberwithin{equation}{section} \theoremstyle{remark}

\newcommand{\ua}{\uparrow}

\title{{\bf Degenerate Fokker-Planck Equations : Bismut Formula, Gradient Estimate  and Harnack Inequality }\footnote{Supported in
 part by NNSFC(11131003), SRFDP, 985 project through the Laboratory of Mathematical and  Complex Systems, and the Fundamental Research Funds for the Central Universities and by ANR EVOL.}
}
\author{
{\bf Arnaud Guillin$^{b)}$ and Feng-Yu Wang$^{a),c)}$\footnote{Corresponding author, wangfy@bnu.edu.cn} }\\
\footnotesize{$^{a)}$ School of Mathematical  Sciences,
Beijing Normal
University, Beijing 100875, China}\\
\footnotesize{$^{b)}$Laboratoire de Math\'ematiques, Universit\'e Blaise Pascal and Institut Universitaire de France, France}\\
 \footnotesize{$^{c)}$ Department of Mathematics,
Swansea University, Singleton Park, SA2 8PP, UK}\\
 \footnotesize{Email: guillin@math.univ-bpclermont.fr;  wangfy@bnu.edu.cn;
F.Y.Wang@swansea.ac.uk}}

\begin{document}
\def\R{\mathbb R}  \def\ff{\frac} \def\ss{\sqrt} \def\B{\scr
B}\def\nn{\nabla} \def\x{\mathbf x} \def\y{\mathbf y} \def\z{\mathbf z}
\def\N{\mathbb N} \def\kk{\kappa} \def\m{{\bf m}}
\def\dd{\delta} \def\DD{\Delta} \def\vv{\varepsilon} \def\rr{\rho}
\def\<{\langle} \def\>{\rangle} \def\GG{\Gamma} \def\gg{\gamma}
  \def\nn{\nabla} \def\pp{\partial} \def\EE{\scr E}
\def\d{\text{\rm{d}}} \def\bb{\beta} \def\aa{\alpha} \def\D{\scr D}
  \def\si{\sigma} \def\ess{\text{\rm{ess}}}
\def\beg{\begin} \def\beq{\begin{equation}}  \def\F{\scr F}
\def\Ric{\text{\rm{Ric}}} \def\Hess{\text{\rm{Hess}}}
\def\e{\text{\rm{e}}} \def\ua{\underline a} \def\OO{\Omega}  \def\oo{\omega}
 \def\tt{\tilde} \def\Ric{\text{\rm{Ric}}}
\def\cut{\text{\rm{cut}}} \def\P{\mathbb P} \def\ifn{I_n(f^{\bigotimes n})}
\def\C{\scr C}      \def\aaa{\mathbf{r}}     \def\r{r}
\def\gap{\text{\rm{gap}}} \def\prr{\pi_{{\bf m},\varrho}}  \def\r{\mathbf r}
\def\Z{\mathbb Z} \def\vrr{\varrho} \def\ll{\lambda}
\def\L{\scr L}\def\Tt{\tt} \def\TT{\tt}\def\II{\mathbb I}
\def\i{{\rm in}}\def\Sect{{\rm Sect}}\def\E{\mathbb E} \def\H{\mathbb H}
\def\M{\scr M}\def\Q{\mathbb Q} \def\texto{\text{o}} \def\LL{\Lambda}
\def\Rank{{\rm Rank}}

\maketitle
\begin{abstract} By constructing successful   couplings for degenerate diffusion processes, explicit derivative formula and   Harnack type inequalities are   presented for   solutions to a class of degenerate Fokker-Planck equations on $\R^m\times\R^{d}$.  The main results are also applied to the study of gradient estimate, entropy/transportation-cost inequality and heat kernel inequalities.\\\smallskip
%\hskip 6.15cm{\bf R\'esum\'e}\\

%Des formules de type Bismut-Elworthy-Li et des in\'egalit\'es de Harnack sont d\'emon\-tr\'ees pour les solutions d'\'equations de Fokker-Planck d\'eg\'en\'er\'ees, ceci par la construction de couplages r\'eussis pour ces diffusions d\'eg\'en\'er\'ees.  Ces r\'esultats sont ensuite appliqu\'es \`a l'\'etude d'estim\'ees (ponctuelles) de gradient des solutions, d'in\'egali\-t\'es entropie-co\^ut de transport ou de noyau de la chaleur.
\end{abstract} \noindent

 \vskip 2cm

 %%%%%%%%%%%%%%%%%%%%%%%%%%%%%%%%%%%%%%%%%%
 %%%%%%%%%%%%%%%%%%%%%%%%%%%%%%%%%%%%%%%%%%

\section{Introduction}

Bismut's derivative formula \cite{Bismut} for  diffusion semigroups on Riemannian manifolds, also known as  Bismut-Elworthy-Li formula due to \cite{EL},  is a powerful tool for stochastic analysis on Riemannian manifolds. On the other hand, the dimension-free Harnack inequality introduced in \cite{W97} has been efficiently applied to the study of functional inequalities, heat kernel estimates and strong Feller properties
in both finite- and infinite-dimensional models, see \cite{ATW06, ATW09, DRW09, ES, K, LW, Ouyang, ORW, RW, W07, W10, WWX10, WX10, WY10, Z}. These two objects have been well developed in the elliptic setting,
but the study for the degenerate case is far from complete.

It is known that the Bismut type formula can be derived for a class of hypoelliptic diffusion semigroups by using Malliavin calculus (see e.g. \cite[Theorem 10]{AT99}). In this case, since no curvature bound can be used, the  derivative formulae are usually less explicit. It is remarkable that in the recent work \cite{Zhang} X. Zhang established an explicit derivative formula for the semigroup associated to degenerate SDEs of type (\ref{1.1}) below (see Section 2 for details). On the other hand, the study of dimension-free Harnack inequality for degenerate diffusion semigroups is very open, except for Ornstein-Uhlenck type semigroups investigated in \cite{ORW}, where the associated stochastic differential equation is linear.

Our strategy is based on coupling, see for example \cite{W10b}, and the main purpose of the paper is thus to construct such a successful coupling using Girsanov transform in the manner of \cite{ATW06} for degenerate diffusion processes, which implies  explicit Bismut formula and dimension-free Harnack inequality for degenerate Fokker-Planck equations.

Let us introduce more precisely the framework we will consider. Let $\si_t$ be  invertible $d\times d$-matrix which is continuous in $t\ge 0$, $A$ be an $m\times d$-matrix with rank $m$,   $B_t$ be a $d$-dimensional Brownian motion, and $Z_t\in C^1(\R^m\times\R^d,\R^d)$ which is continuous in $t$. Consider the following degenerate stochastic differential equation on $\R^m\times \R^d$:

\beq\label{1.1} \beg{cases} \d X_t= A Y_t\d t,\\
\d Y_t= \si_t \d B_t + Z_t(X_t,Y_t)\d t.\end{cases}\end{equation} We  shall use $(X_t(x), Y_t(y))$ to denote the solution
with initial data $(x,y)\in \R^m\times\R^d$.  For simplicity, we will use $\R^{m+d}$ to stand for $\R^{m}\times \R^d$. Then the solution is a Markov process generated by
$$L_t:=\ff 1 2 \sum_{i,j=1}^d (\si_t\si^*_t)_{ij}\ff{\pp^2}{\pp y_i\pp y_j} +\sum_{i=1}^d (Z_t(x,y))_j \ff{\pp}{\pp y_j} +\sum_{l=1}^m (Ay)_l \ff{\pp}{\pp x_l}.$$
For any $f\in \B_b(\R^{m+d}),$   the set of all bounded measurable real functions on $\R^{m+d},$ let
$$ P_t f(x,y):= \E f(X_t(x), Y_t(y)),\ \ t\ge 0, (x,y)\in \R^{m+d}.$$ Then $u(t,x,y):=P_t f(x,y)$ solves the degenerate Fokker-Planck type equation
$$\pp_t u(t,x,y)= L_tu(t,\cdot)(x,y).$$

In the case where $m=d$, $\sigma_t=A=I$ and
$$Z_t(x,y)=-\nabla V(x)-cy,$$
this type of equation has recently attracted much interest under the name $``$kinetic Fokker-Planck equation" in PDE, see Villani \cite{Vil09}, or
 $``$stochastic damping Hamiltonian system" in probability, see \cite{Wu01,BCG08}, where the long time behavior of $P_t$ has been investigated. In this particular case the invariant probability measure (if it exists) is well known as $\mu(\d x,\d y)=\e^{-2V(x)-c|y|^2}\d x \d y$ (up to a constant), and Villani \cite{Vil09} uses this fact to establish hypocoercivity via most importantly an hypoelliptic regularization estimate $H^1\to L^2$.  First note that the methodology used there relies heavily on the knowledge of the invariant measure, which  we will not need in the present study.  Also, his main condition reads as $|\nabla^2V|\le c(1+|\nabla V|)$ preventing exponentially growing potentials, but for parts of our results we do not impose such growing conditions. To allow easier comparison, we will use as running example kinetic Fokker-Planck equation. Let us also mention that we obtain here pointwise estimates, i.e. control of $|\nabla P_t f|$, which allows for example to get uniform bounds when $f$ is initially bounded (exploding when time goes to 0), results that cannot be obtained via Villani's methodology.

 In the following three sections, we will investigate pointwise regularity estimates by establishing derivative formula, gradient estimate  and Harnack inequality for  $P_t$.

 \section{Derivative formulae}

 Since $A$ has rank $m$, we have $d\ge m$ and for any $h_1\in \R^m$, the set
$$A^{-1} h_1:=\{z\in\R^d:\ Az =h_1\}\ne\emptyset.$$
For any $h_1\in\R^m$,  let
$$|A^{-1}h_1|=\inf\{|z|:\ z\in A^{-1}h_1\}.$$ Then it is clear that
$$\|A^{-1}\|:=\sup\big\{|A^{-1}h_1|:\ h_1\in \R^m, |h_1|\le 1\big\}<\infty.$$
We shall use $|\cdot|$ to denote the absolute value and the norm in Euclidean spaces, and use $\|\cdot\|$ to denote
the operator norm of a matrix. For $h\in \R^{m+d}$, we use $D_h$ to stand for the directional derivative along $h$.

Before move on, let us first mention the Bismut formula derived in \cite{Zhang}. We call a $C^2$-function $W$ on $\R^{m+d}$ a Lyapunov function, if $W\ge 1$ having compact level sets. The following result is reorganized  from \cite[Theorem 3.3]{Zhang}. For $h\in \R^{m+d}$, let $\nn_h$ denote the directional derivative along  $h$.

\beg{thm}[\cite{Zhang}] Let $t>0,$ $m=d$ and $A=I$. Assume that there exist a Lyapunov function $W$ and some constants $C>0, \aa\in [0,1], \ll\ge 0$
such that for $s\in [0,t]$
\beq\label{Z0} L_s W\le C W, \ \ |\nn W|^2\le CW^{2-\aa}  \end{equation}  and
\beq\label{Z0'} \beg{cases}|\nn Z_s|\le C W^\ll, \\
\<y-\tt y,Z_s(x,y)-Z_t(\tt x,\tt y)\>\le C|(x-\tt x, y-\tt y)|^2 \big\{W(x,y)^\aa+W(\tt x,\tt y)^\aa\big\} \end{cases}\end{equation} hold for $(x,y),(\tt x,\tt y)\in\R^{m+d}$. Then for any $h=(h_1,h_2)\in \R^{m+d}$ and $f\in \B_b(\R^{m+d})$,
$$\nn_h P_t f= \ff 1 t \E\bigg\{ f(X_t,Y_t)\int_0^t\Big\<\si^{-1}_s\big\{\nn_{\Theta_s} Z_s(X_s,Y_s)-\gg_1'(s)h_1+\gg_2'(s)h_2\big\}, \d B_s\Big\>\bigg\}$$ holds, where
$$\gg_1(s)= 2(t-2s)^++s-t, \ \ \gg_2(s)= \ff 4 t \{s\land (t-s) \}$$ and
$$\Theta_s= \bigg(h_1\int_0^s\gg_1(r)\d r + h_2 t+ h_2 \int_0^s \gg_2(r)\d r,\ \gg_1(s)h_1-\gg_2(s)h_2\bigg).$$\end{thm}

In particular, this result applies to $W(x,y)=1+|x|^2+|y|^2$ and $\aa=0$ provided $|\nn Z|$ is bounded. In general, however,   the assumption $|\nn W|^2\le C W^{2-\aa}$ excludes exponential choices of $W$ like $\exp[|x|^l+|y|^m]$ for $l\lor m>1$, which is exactly the correct Lyapunov function in the study of kinetic Kokker-Planck equation (see
Example 2.1 below). In this section, we aim to present a more general version of the derivative formula without this condition.

Let us introduce now the assumption that we will use in the sequel:

\paragraph{(A)} \emph{There exists a constant $C>0$ such that     $L_s W\le CW$ and
$$|Z_s(\x)-Z_s(\y)|^2\le C|\x-\y|^2 W(\y),\ \ \x,\y\in\R^{m+d}, |\x-\y|\le 1$$
hold for some  Lyapunov function $W$ and $s\in [0,t]$.}

\

Note that condition $L_s W\le CW$, included also in (\ref{Z0}),  is normally a easy to check condition in applications. Although the second condition in {\bf (A)} might be stronger than  (\ref{Z0'}),   it is a natural condition to exchange the order of the expectation and the derivative by using  the dominated convergence theorem, which is however missed in \cite{Zhang} (see line 4 on page 1942 therein).   Most importantly, the second condition in (\ref{Z0}) is now dropped, so that we are able to treat highly non-linear drift $Z$ as in Examples 2.1 and 4.1 below.

The main result in this section provides various different versions of derivative formula by making different choices of the pair functions $(u,v)$.

\beg{thm}\label{T1.1} Assume {\bf (A)}. Then the process $(X_t,Y_t)_{t\ge0}$ is non-explosive for any initial point in $\R^{m+d}$. Moreover, let   $t>0$ and $u,v\in C^2([0,t])$ be such that
\beq\label{1.0} u(t)=v'(0)=1,\ \ \ u(0)=v(0)=u'(0)=u'(t)=v'(t)=v(t)=0.\end{equation} Then for any $h=(h_1,h_2)\in \R^m\times \R^d$ and $z\in A^{-1}h_1:=\{z\in \R^d:\ Az=h_1\},$
\beq\label{1.2} \nn_h P_t f  = \E\bigg\{f(X_t, Y_t)
  \int_0^t \Big\<\si_s^{-1}\big\{u''(s) z -v''(s)h_2
 +(\nn_{\Theta(h,z,s)}Z_s)(X_s,Y_s)\big\},\, \d B_s\Big\>\bigg\}  \end{equation} holds for    $f\in \B_b(\R^{m+d})$,  where
$$\Theta(h,z,s)= \big(\{1-u(s)\}h_1 +v(s)Ah_2,\ v'(s)h_2-u'(s)z\big).$$
\end{thm}
\beg{proof} The non-explosion follows since $L_sW\le CW$ implies
\beq\label{EXP00} \E W(X_s,Y_s)\le W\e^{Cs},\ \ s\in [0,t], (x,y)\in \R^{m+d}.\end{equation}
To prove (\ref{1.2}), we make use of the coupling method with control developed in \cite{ATW06}. Since the process is now degenerate, the construction of coupling is highly technical: we have to force the coupling to be successful before a fixed time by using a lower dimensional  noise.

Let $t>0,(x,y),h=(h_1,h_2)\in \R^{m+d}$ and $z\in A^{-1}h_1$ be fixed. Simply denote   $(X_s,Y_s)=(X_s(x), Y_s(y))$. From now on, let
$$\vv_0= \inf_{s\in [0,t]}\ff 1 {1\lor |\Theta(h,z,s)|}>0,$$ so that   $\vv_0 |\Theta(h,z,s)|\le 1$ for $s\in [0,t]$. For any $\vv\in (0,\vv_0),$
let $(X_s^\vv, Y_s^\vv)$ solve the equation
\beg{equation} \label{CC}\beg{cases} \d  X_s^\vv= A Y_s^\vv\d s,\ \   X_0^\vv= x+\vv h_1,\\
\d  Y_s^\vv =\si_s\d B_s +Z_s(X_s,Y_s)\d s+\vv\{v''(s)h_2-u''(s)z\}\d s,\ \  Y_0^\vv=y+\vv h_2.\end{cases}\end{equation} By (\ref{1.0}) and noting that $Az=h_1$, we have
\beg{equation}\label{2.1}\beg{cases}  Y_s^\vv=Y_s +\vv v'(s)h_2-\vv u'(s)z,\\
 X_s^\vv= x+\vv h_1 + A\int_0^s Y_r^\vv\d r= X_s +\vv \{1-u(s)\}h_1 +\vv v(s)Ah_2.\end{cases}\end{equation}
 Due to (\ref{1.0}), this in particular implies
\beq\label{2.2} (X_t,Y_t)=( X_t^\vv,   Y_t^\vv),\end{equation}
and also that
\beq\label{2.2'}( X_s^\vv,   Y_s^\vv)=(X_s,Y_s)+\vv\Theta(h,z,s),\ \ s\in [0,t].\end{equation}

On the other hand, let
$$\xi^\vv_s= Z(X_s,Y_s)-Z( X_s^\vv, Y_s^\vv) +\vv v''(s) h_2 -\vv u''(s)z,\ \ s\in [0,t]$$ and
\beq\label{WF0}R_s^\vv= \exp\bigg[-\int_0^{s}\<\si^{-1}_s\xi_r^\vv, \d B_r\> -\ff 1 2 \int_0^{s}|\si^{-1}_s\xi_r^\vv|^2\d r\bigg], \ \ s\in [0,t].\end{equation} We have
$$\d Y_s^\vv =\si_s \d B_s^\vv +Z_s( X_s^\vv, Y_s^\vv)\d s$$ for
$$ B_s^\vv:= B_s+ \int_0^s\si^{-1}_s\xi^\vv_r\d r,\ \ s\in [0,t],$$ which is $d$-dimensional Brownian motion under the probability measure
$\Q_\vv:= R^\vv_t\P$ according to Lemma \ref{L1} below and the Girsanov theorem.   Thus, due to (\ref{2.2}) we have
$$ P_t f((x, y)+\vv h) =\E_{\Q_\vv} f( X_t^\vv,  Y_t^\vv) =\E [R^\vv_t f(X_t,Y_t)].$$ Since $P_tf(x,y)=\E f(X_t,Y_t),$ we arrive at
$$ P_tf((x,y)+ \vv h) -P_t f(x,y)= \E[(R^\vv_t-1)f(X_t, Y_t)].$$  The proof is then completed by  Lemma \ref{L2}.
\end{proof}

\beg{lem}\label{L1} If {\bf (A)} holds, then
$$\sup_{s\in [0,t],\vv\in (0,\vv_0)} \E \big(R_s^\vv\log R_s^\vv\big)<\infty.$$Consequently,   for each $\vv\in (0,\vv_0)$,
  $(R_s^\vv)_{s\in [0,t]}$ is a uniformly integrable martingale. \end{lem}

\beg{proof}
$$\tau_n= \inf\{t\ge 0: |X_t(x)|+|Y_t(y)|\ge n\},\ \ n\ge 1.$$ Then $\tau_n\uparrow\infty$ as $n\uparrow\infty.$
By the Girsanov theorem,  $(R_{s\land\tau_n})_{s\in [0,t]}$ is a martingale and $\{B_s^\vv: 0\le s\le t\land \tau_n\}$ is a Brownian motion under the probability measure $\Q_{\vv,n}:=R_{t\land\tau_n}^\vv\P$. Noting that
$$\log R_{s\land \tau_n}^\vv =- \int_0^{s\land\tau_n}\<\si^{-1}_r\xi_r^\vv, \d B_r^\vv\> +\ff 1 2 \int_0^{s\land\tau_n}|\si^{-1}_r\xi_r^\vv|^2\d r, \ \ s\in [0,t],$$
where the stochastic integral is a $\Q_{\vv,n}$-martingale, we have
\beq\label{NNm}
\E[R_{s\land\tau_n}^\vv\log R_{s\land\tau_n}^\vv]=\E_{\Q_{\vv,n}}[\log R_{s\land\tau_n}^\vv]\le
\ff 1 2\E_{\Q_{\vv,n}}   \int_0^{t\land\tau_n}|\si^{-1}_r\xi_r^\vv|^2\d r,\  \ s\in [0,t].\end{equation}
Noting that by {\bf (A)} and (\ref{2.2'})
\beq\label{NNm2} |\si^{-1}_r \xi_r^\vv|^2\le c \vv^2 W(X_r^\vv, Y_r^\vv),\ \ r\in [0,t]\end{equation}  holds for some constant $c>0$, and moreover under the probability measure $\Q_{\vv,n}$
the process $(X_s^\vv, Y_s^\vv)_{s\le t\land\tau_n}$ is generated by $L_s$, $L_sW\le CW$ implies
\beq\label{EXP0} \E_{\Q_{\vv,n} } \int_0^{s\land\tau_n}W(X_r^\vv,Y_r^\vv)\d r \le \int_0^s \E_{\Q_{\vv} }  W(X_r^\vv,Y_r^\vv)\d r\le W(X_0^\vv,Y_0^\vv)\int_0^t\e^{Cr}\d r.  \end{equation}
 Combining this with  (\ref{NNm}) we obtain
\beq\label{WF*}\E[R_{s\land\tau_n}^\vv\log R_{s\land\tau_n}^\vv]\le c,\ \  s\in [0,t], \vv\in (0,\vv_0),n\ge 1\end{equation}    for some constant $c>0.$
Since for each $n$ the process $(R_{s\land\tau_n}^\vv)_{ s\in [0,t]}$ is a martingale, letting $n\to\infty$ in the above inequality  we complete the proof.\end{proof}

\beg{lem}\label{L2} If {\bf (A)} holds then  the family
$\big\{\ff{|R_t^\vv -1|}\vv \big\}_{\vv\in (0,\vv_0)}$  is uniformly integrable w.r.t. $\P$.  Consequently,
\beq\label{FF} \lim_{\vv\to 0} \ff{R_t^\vv-1} \vv = \int_0^t \Big\<\si_s^{-1}\big\{u''(s) z -v''(s)h_2  +(\nn_{\Theta(h,z,s)}Z)(X_s(x),Y_s(y))\big\},\, \d B_s\Big\> \end{equation}  holds
 in $L^1(\P).$
\end{lem}

\beg{proof}  Let $\tau_n$ be in the proof of Lemma \ref{L1} and let
$$N_s^\vv= \si_s^{-1} \big\{\nn_{\Theta(h,z,s)}Z_s(X_s^\vv,Y_s^\vv)+u''(s)z-v''(s)h_2\big\},\ \ s\in [0,t], \vv\in (0,\vv_0).$$ By {\bf (A)} and (\ref{NNm2}), there exists a constant $c>0$ such that
\beq\label{WFF} \big|\< N_s^\vv, \si_s^{-1} \xi_s^\vv\>\big|\le \vv |N_s^\vv|^2+\vv^{-1}|\si^{-1}_s\xi_s^\vv|^2\le c\vv W(X_s^\vv,Y_s^\vv),\ \
\vv\in (0,\vv_0), s\in [0,t].\end{equation}
 Since $\nn Z$ is locally bounded, it follows from (\ref{2.2'}) and (\ref{WF0}) that
$$ \ff{\d}{\d\vv} R_{t\land \tau_n}^\vv  = R_{t\land\tau_n}^\vv\bigg\{\int_0^{t\land\tau_n}\<N_s^\vv, \d B_s\> +  \int_0^{t\land \tau_n}\<N_s^\vv, \si_s^{-1} \xi_s^\vv\>\d s\bigg\},\ \ \vv\in (0,\vv_0), n\ge 1.$$
Combining this with (\ref{WFF}) we obtain
$$  \ff{|R_{t\land\tau_n}^\vv -1|}\vv   \le \ff 1 \vv \int_0^\vv R_{t\land\tau_n}^r\d r  \int_0^{t\land\tau_n}\<N_s^r, \d B_s\>
 + c\int_0^{\vv_0} R_{t\land \tau_n}^r\d r \int_0^{t\land\tau_n} W(X_s^r,Y_s^r)\d s$$ for $\vv\in (0,\vv_0),n\ge 1.$
Noting that
under $\Q_{r}$ the process $(X_s^r,Y_s^r)_{s\in [0,t]}$ is generated by $L_s$, by (\ref{EXP00}) we have
$$\E \int_0^{\vv_0} R_{t}^r\d r \int_0^{t}W(X_s^r,Y_s^r)\d s= \int_0^{\vv_0} \d r \int_0^t \E_{\Q_r} W(X_s^r,Y_s^r)\d s <\infty.$$ Thus, for the first assertion it remains to show that the family
$$\eta_{\vv,n}:= \ff 1 \vv \int_0^\vv R_{t\land \tau_n}^r |\Xi_{t,n}|(r)\d r,\ \ \vv\in (0,\vv_0), n\ge 1$$ is uniformly integrable, where
$$\Xi_{t,n}(r):=  \int_0^{t\land\tau_n}\<N_s^r, \d B_s\>.$$
Since $r\log^{1/2} (\e+r)$ is increasing and convex in $r\ge 0$, by the Jensen inequality,
\beg{equation*}\beg{split} &\E\big\{\eta_{\vv,n} \log^{1/2}(\e+\eta_{\vv,n})\big\} \\
&\le \ff 1\vv \int_0^\vv \E\Big\{R_{t\land\tau_n}^r |\Xi_{t,n}|(r)\log^{1/2}\big(\e + R_{t\land\tau_n}^r |\Xi_{t,n}|(r)\big)\Big\}\d r\\
&\le \ff 1 \vv \int_0^\vv \E\Big\{R_{t\land\tau_n}^r |\Xi_{t,n}|(r)^2 + R_{t\land\tau_n}^r  \log\big(\e + R_{t\land\tau_n}^r |\Xi_{t,n}|(r)\big)\Big\}\d r\\
&\le \ff 1\vv\int_0^\vv \E\Big\{c+ 2 R_{t\land\tau_n}^r |\Xi_{t,n}|(r)^2 +  R_{t\land\tau_n}^r\log  R_{t\land\tau_n}^r\Big\}\d r\end{split}\end{equation*} holds for some constant $c>0$. Combining this with   (\ref{WF*}) and noting that (\ref{WFF}) and (\ref{EXP0}) imply
\beg{equation*}\beg{split} &\E \big\{R_{t\land\tau_n}^r |\Xi_{t,n}|(r)^2\big\}= \E_{\Q_{r,n}}\bigg( \int_0^{t\land\tau_n}
\big\<N_s^r,\ \d B_s^r\big\>\bigg)^2
 =\E_{\Q_{r,n}}\int_0^{t\land\tau_n}|N_s^r|^2\d s \\
 &\le c \E_{\Q_{r,n}}\int_0^{t\land \tau_n}  W(X_s^r, Y_s^r)\d s\le c',\ \ \ n\ge 1, r\in(0,\vv_0) \end{split}\end{equation*} for some constants $c,c'>0$,    we conclude that $\{\eta_{\vv,n}\}_{\vv\in (0,\vv_0),n\ge 1}$ is uniformly integrable. Thus, the proof of the first assertion is finished.

Next, by {\bf (A)} and (\ref{2.2'})  we have
$$\lim_{\vv\to 0} \Big|\ff {\xi_s^\vv}\vv + (\nn_{\Theta(h,z,s)}Z)(X_s,Y_s)+u''(s)z-v''(s)h_2\Big|=0.$$ Moreover, for each $n\ge 1$ this sequence
is bounded on $\{\tau_n\ge t\}$. Thus, $(\ref{FF})$ holds a.s. on $\{\tau_n\ge t\}$. Since $\tau_n\uparrow\infty$, we conclude that (\ref{FF}) holds a.s.
Therefore, it also holds on $L^1(\P)$ since  $\{\ff{R_t^\vv-1}\vv\}_{\vv\in (0,1)}$ is uniformly integrable according to the first assertion. \end{proof}

To conclude this section, we present an example of kinetic Fokker-Planck equation for which $W$ is an exponential function so that (\ref{Z0'}) fails true but {\bf (A)} is satisfied.

\paragraph{Example 2.1}(Kinetic Fokker-Planck equation)
Let $m=d $ and consider
\beq\label{kfp} \beg{cases} \d X_t= Y_t\d t,\\
\d Y_t= \d B_t -\nabla V(X_t)\d t-Y_t\d t\end{cases}\end{equation}
for some $C^2$-function $V\ge 0$ with compact let sets. Let $W(x,y)=\exp[2V(x)+|y|^2]$. We easily get that
$L W= d  W.$ Thus, it is easy to see that {\bf (A)} holds for e.g. $V(x)= (1+|x|^2)^l$  or even $V(x)=\e^{(1+|x|^2)^l}$ for some constant $l\ge 0$. Therefore, by Theorem \ref{T1.1}  the derivative formula (\ref{1.2}) holds for $(u,v)$ satisfying (\ref{1.0}).\\
Note that Villani \cite[th. A.8]{Vil09} has a crucial assumption: $|\nabla^2 V|\le C(1+|\nabla V|)$ which prevents potential behaving as $V(x)=\e^{(1+|x|^2)^l}$. Note also that the previous arguments do not rely on the explicit knowledge of an invariant probability measure, which is crucial in Villani's argument.

\section{Gradient estimates}

In this section we aim to derive gradient estimates from the derivative formula   (\ref{1.2}). For simplicity, we only consider the time-homogenous case that $\si$ and $Z$ are independent of $t$. In general, we have the following result.

\beg{prp} \label{P3.1} Assume {\bf (A)} and let $(u,v)$ satisfy $(\ref{1.0})$. Then for any $f\in \B_b(\R^{m+d}), t>0$ and $h=(h_1,h_2)\in\R^{m+d}$,
$z\in A^{-1} h_1$,
\beq\label{3.1} |\nn_h P_t f|^2 \le \|\si^{-1}\|^2(P_t f^2) \E \int_0^t \big|u''(s)z-v''(s)h_2+\nn_{\Theta(h,z,s)}Z(X_s,Y_s)\big|^2\d s.\end{equation}If $f\ge 0$ then
for any $\dd>0$,
\beq\label{3.2}\beg{split}& |\nn_h P_t f|\le  \dd\big\{P_t (f\log f)-(P_t f)\log P_t f\big\}\\
&\qquad+\ff{\dd P_t f}2   \log \E\exp\bigg[\ff {2\|\si^{-1}\|^2} {\dd^2} \int_0^t\big|u''(s)z-v''(s)h_2+\nn_{\Theta(h,z,s)}Z(X_s,Y_s)\big|^2\d s\bigg].\end{split}\end{equation}
\end{prp}
\beg{proof} Let $M_t= \int_0^t \big\<\si^{-1} \big\{u''(s)z-v''(s)h_2+\nn_{\Theta(h,z,s)}Z(X_s,Y_s)\big\},\ \d B_s\big\>.$  By
  (\ref{1.2}) and the Schwartz inequality we obtain
$$ |\nn_h P_t f|^2 \le  (P_t f^2) \E M_t^2  \le \|\si^{-1}\|^2(P_t f^2) \E \int_0^t \big|u''(s)z-v''(s)h_2+\nn_{\Theta(h,z,s)}Z(X_s,Y_s)\big|^2\d s.$$ That is, (\ref{3.1}) holds. Similarly, (\ref{3.2}) follows from (\ref{1.2}) and the Young inequality (cf. \cite[Lemma 2.4]{ATW09}):
$$|\nn_h P_t f| \le \dd \big\{P_t (f\log f)-(P_t f)\log P_t f\big\}+\dd \log\E\exp\Big[\ff{M_t}\dd\Big]$$ since
\beg{equation*}\beg{split} &\E\exp\Big[\ff{M_t}\dd\Big]\le \bigg(\E\exp\Big[\ff{2\<M\>_t}{\dd^2}\Big]\bigg)^{1/2}\\
&\le \bigg(\E\exp\bigg[\ff {2\|\si^{-1}\|^2} {\dd^2} \int_0^t\big|u''(s)z-v''(s)h_2+\nn_{\Theta(h,z,s)}Z(X_s,Y_s)\big|^2\d s\bigg]\bigg)^{1/2}.\end{split}\end{equation*}\end{proof}

To derive explicit estimates, we will take the following explicit choice of the
 pair $(u, v)$:
\beq\label{uv1}u(s)= \ff{s^2(3t-2s)}{t^3},\ \ \ v(s)=  \ff{s(t-s)^2}{t^2},\ \ \ \ s\in [0,t],\end{equation}  which satisfies  (\ref{1.0}). In this case
we have
\beg{equation}\label{uv2}\beg{split} &u'(s)=\ff{6s(t-s)}{t^3},\ u''(s)= \ff{6(t-2s)}{t^3},  \ v'(s)=\ff{(t-s)(t-3s)}{t^2},\\
&\ v''(s)=\ff{2(3s-2t)}{t^2},\ 1-u(s)= \ff{(t-s)^2(t+2s)}{t^3} ,\ \ s\in [0,t].\end{split}\end{equation}
In this case, Proposition \ref{P3.1} holds for
\beq\label{uv}\beg{split} &u''(s)z-v''(s)h_2=\LL(h,z,s):= \ff{6(t-2s)}{t^3}z +\ff{2(2t-3s)}{t^2}h_2,\\
&\Theta(h,z,s)= \\
&\ \Big(\ff{(t-s)^2(t+2s)}{t^3}h_1+\ff{s(t-s)^2}{t^2} A h_2, \ff{(t-s)(t-3s)}{t^2}h_2-\ff{6s(t-s)}{t^3}z\Big).\end{split}\end{equation}

Below we consider the following three cases respectively:
\beg{enumerate}\item[(i)] $|\nn Z|$ is bounded;
\item[(ii)] $|\nn Z|$ has polynomial growth and $\<Z(x,y), y\>\le C(1+|x|^2+|y|^2)$ holds for some constant $c>0;$
\item[(iii)] A more general case including the kinetic Fokker-Planck equation.\end{enumerate}

\subsection{Case (i): $|\nn Z|$ is bounded}

In this case {\bf (A)} holds for e.g. $W(x,y)=1+|x|^2+|y|^2$, so that Proposition \ref{P3.1} holds for $u''(s)z-v''(s)h_2$ and $\Theta(h,z,s)$ given in (\ref{uv}). From this specific choice of $\Theta(h,z,s)$
we see that $\nn^x Z$ and $\nn^y Z$   will lead  to different time behaviors of   $\nn_h P_t f$.  So, we adopt the condition
\beq\label{NZ} |\nn^x Z(x,y)|\le K_1,\ \ |\nn^y Z(x,y)|\le K_2,\ \ (x,y)\in \R^{m+d}\end{equation} for some constants $K_1,K_2\ge 0,$ where $\nn^x$ and $\nn^y$ are the gradient operators w.r.t. $x\in\R^m$ and $y\in\R^d$ respectively.
Moreover, for $t>0$ and $r_1,r_2\ge 0$, let
$$\Psi_t(r_1,r_2)=\|\si^{-1}\|^2t\bigg\{r_1\Big(\ff {6\|A^{-1}\|}{t^2} +K_1 +\ff{3K_2\|A^{-1}\|}{2t}\Big)
+r_2 \Big(\ff 4 t +\ff{4K_1 t\|A\|}{27}+K_2\Big)\bigg\}^2$$ and
\beq\label{Phi} \Phi_t(r_1,r_2)= \inf_{s\in (0,t]} \Psi_s(r_1,r_2).\end{equation}
In the following result  the inequality (\ref{1.3}) corresponds to the pointwise estimate of the $H^1\to L^2$ regularization investigated in Villani \cite[Th. A.8]{Vil09}, while (\ref{1.5}) corresponds to the pointwise estimate of the regularization $``$Fisher information to entropy" \cite[Th A.18]{Vil09}.

\beg{cor} \label{C2} Let $(\ref{NZ})$ hold for some constants $K_1,K_2\ge 0$. Then for any $t>0, h=(h_1,h_2)\in\R^{m+d}$,
\beq\label{1.3} |\nn_h P_t f|^2 \le   (P_t f^2)\Phi_t(h_1,h_2),\ \  f\in \B_b(\R^{m+d}).
\end{equation} If  $f\ge 0$, then
\beq\label{1.4}|\nn_h P_t f|\le \dd \big\{P_t (f\log f)- (P_t f)\log (P_tf)\big\}
 +\ff {P_tf}\dd \Phi_t(h_1,h_2)\end{equation} holds for all $\dd>0,$ and consequently
  \beq\label{1.5}|\nn_h P_t f|^2\le 4\Phi_t(h_1,h_2) \big\{P_t (f\log f)- (P_t f)\log (P_tf)\big\} P_t f. \end{equation}
\end{cor}

\beg{proof}Let $z$ be such that $|z|=|A^{-1}h_1|\le \|A^{-1}\||h_1|,$ and take
\beq\label{Eta} \eta_s= \LL(h,z,s)  +\nn_{\Theta(h,z,s)}Z(X_s(x),Y_s(y)).\end{equation}
By (\ref{3.1}),
\beq\label{3.1'} |\nn_h P_tf (x,y)|^2\le \|\si^{-1}\|^2(P_tf^2)(x,y)\E\int_0^t|\eta_s|^2\d s.\end{equation}
Since (\ref{NZ}) implies $|\nn_hZ|\le K_1|h_1|+K_2|h_2|$,  it follows  that
\beg{equation*}\beg{split} |\eta_s|&\le \Big|\ff{6(t-2s)}{t^3}z +\ff{2(2t-3s)}{t^2}h_2\Big| +K_1\Big|\ff{(t-s)^2(t+2s)}{t^3}h_1+\ff{s(t-s)^2}{t^2} A h_2\Big|\\
&\qquad\qquad
+K_2\Big|\ff{(t-s)(t-3s)}{t^2}h_2-\ff{6s(t-s)}{t^3}z\Big|\\
&\le |h_1|\Big(\ff {6\|A^{-1}\|}{t^2} +K_1 +\ff{3K_2\|A^{-1}\|}{2t}\Big)
+|h_2| \Big(\ff 4 t +\ff{4K_1 t\|A\|}{27}+K_2\Big).\end{split}\end{equation*} Then
\beq\label{3.2'}\beg{split} \int_0^t | \eta_s|^2\d s\le&
   t \bigg\{|h_1|\Big(\ff {6\|A^{-1}\|}{t^2}+K_1 +\ff{3K_2\|A^{-1}\|}{2t}\Big)\\&\qquad\qquad\qquad
+|h_2| \Big(\ff 4 t +\ff{4K_1 t\|A\|}{27}+K_2\Big)\bigg\}^2.\end{split}\end{equation}
Combining this with (\ref{3.1'}) we obtain
$$  |\nn_h P_t f|^2 \le   (P_t f^2)\Psi_t(|h_1|,|h_2|).$$
Therefore,  for any $s\in (0,t]$ by the semigroup property and the Jensen inequality one has
$$|\nn P_t f|^2=|\nn P_s(P_{t-s} f)|^2\le \Psi_s(|h_1|,|h_2|)P_s(P_{t-s}f)^2\le \Psi_s(|h_1|,|h_2|)P_tf^2.$$
This proves (\ref{1.3}) according to  (\ref{Phi}).

To prove (\ref{1.4}) we let $f\ge 0$ be bounded. By (\ref{3.2}),
\beq\label{EG}\beg{split} |\nn_hP_tf|\le &\dd \big\{P_t(f\log f)-(P_tf)\log (P_tf)\big\}\\
 &+\ff {\dd P_tf} 2 \log \E\,\exp\bigg[\ff {2\|\si^{-1}\|^2} {\dd^2} \int_0^t |\eta_s|^2\d s\bigg].\end{split}\end{equation}
Combining this with (\ref{3.2'}) we obtain
$$ |\nn_hP_tf|\le   \dd \big\{P_t(f\log f)-(P_tf)\log (P_tf)\big\}
  +\ff {P_tf}\dd \Psi_t(|h_1|,|h_2|).$$  As observed above,  by the semigroup property and the Jensen inequality, this  implies (\ref{1.4}).

 Finally, minimizing the right hand side of (\ref{1.4}) in $\dd>0$, we obtain
 $$|\nn_h P_t f|\le 2 \ss{\Phi_t (|h_1|,|h_2|)\{P_t(f\log f)-(P_t f)\log P_t f\}P_t f}.$$ This is equivalent to (\ref{1.5}). \end{proof}

\subsection{Case (ii)}
Assume there exists $l> 0$ such that

\paragraph{(H)} (i)  $\<Z(x,y), y\> \le C(|x|^2+|y|^2+1),\ (x,y)\in \R^{m+d}$;

\ \ \ (ii) $|\nn Z|(x,y):=\sup\{|\nn_h Z|(x,y):\ |h|\le 1\}\le C(1+|x|^2+|y|^2)^l,\ (x,y)\in\R^{m+d}.$

\

It is easy to see that {\bf (H)} implies {\bf (A)} for $W(x,y)=(1+|x|^2+|y|^2)^{2l},$ so that Proposition \ref{P3.1} holds for $u''(s)z-v''(s)h_2$ and $\Theta(h,z,s)$ given in (\ref{uv}).

\beg{cor} \label{C3} Let {\bf (H)} hold. \beg{enumerate} \item[$(1)$] There exists a constant $c>0$ such that
$$ |\nn P_tf|^2 (x,y)\le \ff{c}{(t\land 1)^3} P_t f^2(x,y),\ \ f\in \B_b(\R^{m+d}),\ \ t>0, (x,y)\in \R^{m+d}.$$
\item[$(2)$] If   $l<\ff 1 2$, then there exists a constant $c>0$ such that
\beg{equation*}\beg{split}  |\nn P_t f|(x,y)\le &\dd \big\{P_t (f\log f)- (P_t f)\log (P_tf)\big\}(x,y)\\
&+ \ff {c P_t f(x,y)} {\dd (t\land 1)^4}\big\{(|x|^2+|y|^2)^{2l} + (\dd(1\land  t)^2)^{4(l-1)/(1-2l)}\big\}\end{split}\end{equation*} holds for all $\dd>0$ and positive $f\in \B_b(\R^{m+d})$ and $(x,y)\in\R^{m+d}$.
\item[$(3)$] If $l=\ff 1 2$, then there exist two constants $c,c'>0$ such that for any $t>0$ and $\dd\ge t^{-2}\e^{c(1+t)}$,
$$ |\nn P_t f|(x,y)\le  \dd\big\{P_t(f\log f)-(P_t f)\log P_t f\big\}(x,y) + \ff{c'  P_tf(x,y)}\dd \big( 1+|x|^2+|y|^2\big)$$ holds for all positive $f\in \B_b(\R^{m+d}) $ and $(x,y)\in \R^{m+d}.$
\end{enumerate}
\end{cor}
\beg{proof}   As observed in the proof of Corollary \ref{C2}, we only have to prove the results for $t\in (0,1].$

(1) It is easy to see that $\eta_s$ in the proof of Corollary \ref{C2} satisfies
\beq\label{si}|\si^{-1}\eta_s|^2\le c_1(t^2+t^{-4}) |h|^2 (1+|X_s(x)|^2+|Y_s(y)|^2)^{2l}\end{equation}
for some constant $c_1>0$. Thus, the first assertion follows from (\ref{3.1'}) and Lemma \ref{L1}.

(2)  Let {\bf (H)} hold for some $l\in (0,1/2)$. Then $$L(1+|x|^2+|y|^2)^{2l}\le c_2 (1+|x|^2+|y|^2)^{2l}$$ holds for some constant $c_2>0$.
Let $(X_s, Y_s)= (X_s(x), Y_s(y))$. By the It\^o formula, we have
$$\d (1+|X_s|^2+|Y_s|^2)^{2l} \le 4l  (1+|X_s|^2+|Y_s|^2)^{2l-1}\<Y_s,\si\d B_s\> +c_2  (1+|X_s|^2+|Y_s|^2)^{2l}\d s.$$ Thus,
\beg{equation*}\beg{split} &\d\big\{\e^{-(1+c_2)s} (1+|X_s|^2+|Y_s|^2)^{2l}\big\}\\
& \le 4l \e^{-(1+c_2)s} (1+|X_s|^2+|Y_s|^2)^{2l-1}\<Y_s,\si\d B_s\> -
\e^{-(1+c_2)s} (1+|X_s|^2+|Y_s|^2)^{2l}\d s.\end{split}\end{equation*} Therefore, for any $\ll>0$,
\beg{equation}\label{EXP}\beg{split} &\E \e^{\ll \int_0^t \e^{-(1+c_2)s}  (1+|X_s|^2+|Y_s|^2)^{2l}\d s}\\
& \le \e^{\ll (1+|x|^2+|y|^2)^{2l}}
\E \e^{4\ll l\int_0^t \e^{-(1+c_2)s}  (1+|X_s|^2+|Y_s|^2)^{2l-1}\<Y_s,\si \d B_s\>}\\
&\le \e^{\ll (1+|x|^2+|y|^2)^{2l}}\Big\{\E \e^{16\ll^2l^2\|\si\|^2 \int_0^t \e^{-2(1+c_2)s}  (1+|X_s|^2+|Y_s|^2)^{2(2l-1)}|Y_s|^2\d s}\Big\}^{1/2}\\
&\le \e^{\ll (1+|x|^2+|y|^2)^{2l}}\Big\{\E \e^{16\ll^2l^2\|\si\|^2 \int_0^t \e^{-(1+c_2)s}  (1+|X_s|^2+|Y_s|^2)^{4l-1}\d s}\Big\}^{1/2}.\end{split}\end{equation}
On the other hand, since $l<\ff 1 2$ implies $4l-1< 2l$, there exists a constant $c_3>0$ such that
$$16\ll^2l^2\|\si\|^2r^{4l-1}\le \ll r^{2l}+ c_3 \ll^{(3-4l)/(1-2l)},\ \ \ r\ge 0.$$ Combining this with (\ref{EXP}) we arrive at
\beg{equation*}\beg{split}   & \E \exp\bigg[\ll \int_0^t \e^{-(1+c_2)s}  (1+|X_s|^2+|Y_s|^2)^{2l}\d s\bigg]\\
& \le  \exp\Big[\ll (1+|x|^2+|y|^2)^{2l}+\ff{c_3}2\ll^{(3-4l)/(1-2l)}\Big]\\
&\quad\times \bigg(  \E \exp \bigg[ \ll \int_0^t \e^{-(1+c_2)s}  (1+|X_s|^2+|Y_s|^2)^{2l}\d s\bigg]\bigg)^{1/2}.\end{split}\end{equation*}
As the argument works also for $t\land \tau_n$ in place of  $t$, we may   assume priorly that the left-hand side of the above inequality is finite, so that
\beg{equation*} \E \exp\bigg[\ll \int_0^t \e^{-(1+c_2)s}  (1+|X_s|^2+|Y_s|^2)^{2l}\d s\bigg]\le \exp\Big[2\ll (1+|x|^2+|y|^2)^{2l}+c_3\ll^{(3-4l)/(1-2l)}\Big].\end{equation*}
Letting $$\ll_t(\dd)= \ff{2c_1(t^2+t^{-4})}{\dd^2} \e^{(1+c_2)t},$$ and combining the above inequality with (\ref{EG}) and (\ref{si}), we arrive at
\beg{equation}\label{EXP2}\beg{split} &\Big(|\nn P_t f|-\dd\big\{P_t(f\log f)-(P_t f)\log P_t f\big\}\Big)(x,y)\\
&\qquad\le \ff {\dd P_t f(x,y)}2\log \E\exp\bigg[\ll_t(\dd) \int_0^t \e^{-(1+c_2)s}  (1+|X_s|^2+|Y_s|^2)^{2l}\d s\bigg]\\
&\qquad\le \dd P_t f(x,y)\Big\{ \ll_t(\dd) (1+|x|^2+|y|^2)^{2l} +\ff {c_3} 2 \ll_t(\dd)^{(3-4l)/(1-2l)}\Big\}\\
&\qquad\le  \ff {P_t f(x,y) \e^{c(1+t)}} {\dd t^4}\big\{(|x|^2+|y|^2)^{2l} + \dd^{4(l-1)/(1-2l)}t^{8(l-1)/(1-2l)}\big\}\end{split}\end{equation}
for some constant $c>0$. This proves the desired estimate for $t\in (0,1]$, and hence for all $t>0$ as observed in the proof of Corollary \ref{C2}.

(3) Let {\bf (H)} hold for $l=\ff 1 2$, so that (\ref{EXP}) reduces to
$$\E \e^{\ll \int_0^t \e^{-(1+c_2)s}  (1+|X_s|^2+|Y_s|^2)\d s}\le
\e^{\ll (1+|x|^2+|y|^2)}\Big\{\E \e^{4\ll^2\|\si\|^2 \int_0^t \e^{-(1+c_2)s}  (1+|X_s|^2+|Y_s|^2)\d s}\Big\}^{1/2}.$$
Taking $\ll= (2\|\si\|)^{-2}$ we obtain
$$\E  \exp\bigg[\ff 1 {4\|\si\|^2} \int_0^t \e^{-(1+c_2)s}  (1+|X_s|^2+|Y_s|^2)\d s\bigg]\le
\exp\Big[\ff 1 {4\|\si\|^2} (1+|x|^2+|y|^2)\Big].$$ Obviously, there exists a constant $c>0$ such that if $\dd\ge t^{-2}\e^{c(1+t)}$ then
$\ll_t(\dd)\le (2\|\si\|)^{-2}$ so that

\beg{equation*}\beg{split} &\Big(|\nn P_t f|-\dd\big\{P_t(f\log f)-(P_t f)\log P_t f\big\}\Big)(x,y)\\
&\qquad\le  \ff {\dd P_t f(x,y)}2\log \bigg(\E\exp\bigg[\ff 1 {4\|\si\|^2} \int_0^t \e^{-(1+c_2)s}  (1+|X_s|^2+|Y_s|^2)^{2l}\d s\bigg]\bigg)^{4\|\si\|^2/\ll_t(\dd)}\\
&\qquad\le \ff {\dd P_t f(x,y)}{2\ll_t(\dd)}\big(1+|x|^2+|y|^2\big)\le\ff{c'  P_tf(x,y)}\dd \big( 1+|x|^2+|y|^2\big)\end{split}\end{equation*}holds for some constant $c'>0$.
\end{proof}

\paragraph{Example 3.1} (Kinetic Fokker-Planck equation)

 Let us consider once again the Example 2.1 introduced previously, and remark that the result of Corollary \ref{C3} (1) holds without the first assumption in {\bf (H)}, so that we get a pointwise version of Villani \cite[Th. A.8]{Vil09} under the same type of condition (polynomial growth at most), and thus recover its $L_2$ bound (constants are however rather difficult to compare).

\subsection{A  general case}
\beg{cor}\label{C4} Assume {\bf (A)}. Then there exists a constant $c>0$ such that
\beq\label{A1} |\nn P_t f|^2 \le c\Big(\ff 1 {(1\land t)^3}+ \ff W{1\land t}\Big) P_t f^2,\ \  f\in \B_b(\R^{m+d}).\end{equation}
If moreover there exist constants $\ll, K>0$ and a $C^2$-function $\tt W\ge 1$ such that
\beq\label{A2} \ll W\le K-\ff{L\tt W}{\tt W},\end{equation} then there exist constants $c,\dd_0>0$ such that
\beq\label{GGA} |\nn P_t f|\le \dd \big\{P_t(f\log f)-(P_t f)\log P_t f\big\} + \ff{c}{\dd}\Big\{ \ff 1 {(t\land 1)^3} +\ff{\log \tt W}{(t\land 1)^2}\Big\} P_t f\end{equation} holds for
$f\in \B_b^+(\R^{m+d})$ and $\dd\ge \dd_0/t$.  \end{cor}

\beg{proof} Again, it suffices to prove for $t\in (0,1].$ By (\ref{uv}) and taking $z\in A^{-1}h_1$ such that $|z|=\|A^{-1}\|\cdot |h_1|$, there exists a constant $c>0$ such that
$$|\LL(h,z,s)|\le \ff c {t^2} |h|,\ \ |\Theta(h,z,s)|\le \ff c t|h|.$$ So, by {\bf (A)}
\beq\label{AA0} \big|u''(s)z-v''(s)h_2+\nn_{\Theta(h,z,s)}Z(X_s,Y_s)\big|^2\le \ff c {t^4} +  \ff c {t^2} W(X_s,Y_s)\end{equation} holds for some constant $c>0.$ Since $W\ge 1$ and  $\E W(X_s,Y_s)\le \e^{Cs}W$,
 this and  (\ref{3.1}) yield that
\beg{equation*}\beg{split} |\nn P_t f|^2  &\le c_1 (P_t f^2)\bigg\{\int_0^t |\LL(h,z,s)|^2\d s + \E\int_0^t |\Theta(h,z,s)|^2W(X_s,Y_s)\d s\bigg\}\\
 &\le c_2\Big(\ff 1 {t^3} +\ff W t\Big)  P_t f^2 \end{split}\end{equation*}  holds for some constants $c_1,c_2>0.$

Next,  it is easy to see that  the process
$$M_s:=\tt W(X_s,Y_s) \exp\bigg[-\int_0^s \ff{L\tt W}{\tt W}(X_r,Y_r)\d r\bigg]$$ is a local martingale, and thus a supermartingale due to the Fatou lemma. Combining this with (\ref{A2}) and noting that $\tt W\ge 1$, we obtain
\beq\label{*AA}\E\e^{\ll \int_0^t W(X_s, Y_s)\d s} \le \e^{ K t}\E M_t\le \e^{K t} \tt W.\end{equation} Then the second assertion follows from (\ref{3.2}) and (\ref{AA0}) since for any constant $\aa>0$   there exists a constant $c_2>0$ such that for any   $\dd t\ge \ss{\aa/\ll}$,
$$\E\exp\bigg[\ff{\aa}{\dd^2t^2} \int_0^tW(X_s,Y_s)\d s\bigg]
 \le  \bigg(\E \exp\bigg[\ll \int_0^t W(X_s,Y_s)\d s \bigg]\bigg)^{\aa/(\ll\dd^2t^2)}. $$
 \end{proof}

%%%%%%%%%%%%%%

\section{Harnack inequality and applications}

The aim of this section is to establish the log-Harnack inequality introduced in \cite{RW, W10} and the Harnack  inequality with power due to \cite{W97}. Applications of these inequalities to heat kernel estimates as well as Entropy-cost inequalities can be found in e.g. \cite{RW,W10}. We first consider the general case with assumption {\bf (A)} then move to the more specific setting with assumption {\bf (H)}. Again, we only consider
the time-homogenous case.

\subsection{Harnack inequality under {\bf (A)}}

We first introduce a result, essentially due to \cite{ATW09}, that the entropy-gradient estimate (\ref{3.2}) implies the Harnack inequality with a power.

\beg{prp}\label{P4.0} Let $\H$ be a Hilbert space and $P$ a Markov operator on $\B_b(\H).$ Let $h\in\H$ such that for some   $\dd_h\in (0,1)$ and measurable function $\gg_h: [\dd_h,\infty)\times \H\to (0,\infty)$,
\beq\label{EGF} |\nn_h Pf|\le \dd \big\{P(f\log f)- (Pf)\log Pf\big\} +\gg_h(\dd,\cdot)Pf,\ \ \dd\ge \dd_h\end{equation}holds for all positive $f\in\B_b(\H)$.
Then for any $\aa\ge \ff 1 {1-\dd_h}$ and positive $f\in\B_b(\H)$,
$$(Pf)^\aa(\x)\le  Pf^\aa(\x+h)  \exp\bigg[\int_0^1 \ff{\aa}{1+(\aa-1)s}\gg_h\Big(\ff{\aa-1}{1+(\aa-1)s}, \x+sh\Big)\d s\bigg],\ \ \x\in\H.$$\end{prp}

\beg{proof} Let $\bb(s)= 1+(\aa-1)s.$ We have $\ff{\aa-1}{\bb(s)}\ge \dd_h$ provided $\aa\ge \ff 1 {1-\dd_h}.$ Then
\beg{equation*}\beg{split} &\ff{\d}{\d s} \log (P f^{\bb(s)})^{\aa/\bb(s)}(\x+sh)\\
&=\ff{\aa(\aa-1)\{P (f^{\bb(s)}\log f^{\bb(s)})-(P f^{\bb(s)})\log P f^{\bb(s)}\}}{\bb(s)^2P f^{\bb(s)}}(\x+sh)
 +\ff{\aa \nn_h P f^{\bb(s)}}{\bb(s)P f^{\bb(s)}}(\x+sh)\\
&\ge -\ff{\aa}{\bb(s)} \gg_h\Big(\ff{\aa-1}{\bb(s)},\ \x+sh\Big),\ \ \ s\in [0,1].\end{split}\end{equation*} Then the proof is completed by taking integral over $[0,1]$ w.r.t. $\d s$. \end{proof}

Below is a  consequence of (\ref{GGA}) and Proposition \ref{P4.0}.

\beg{cor}\label{C0} Let {\bf (A)} and $(\ref{A2})$ hold. Then there exist constants $\dd_0,c>0$ such that for any  $\aa>1, t>0$ and positive  $f\in \B_b(\R^{m+d})$,

\beq\label{HHH} (P_t f)^\aa(\x)\le P_t f^\aa(\x+h) \exp\bigg[ \ff{\aa c |h|^2}{\aa-1}\bigg(\ff 1 {(1\land t)^3}+\ff {\int_0^1 \log\tt W(\x+sh)\d s} {(1\land t)^2}\bigg)\bigg]\end{equation} holds for $\x,h\in \R^{m+d}$ with $|h|<\dd_0 t.$\end{cor}

\beg{proof}   By (\ref{GGA}),
$$|\nn_h P_t f|\le \dd |h|\big\{P_t(f\log f)-(P_t f)\log P_t f\big\}
+ \ff{c}{\dd}\Big\{ \ff 1 {(t\land 1)^3} +\ff{\log \tt W}{(t\land 1)^2}\Big\} P_t f$$ holds for $\dd\ge \dd_0/t$. Thus, (\ref{EGF}) holds for $P=P_t$ and
$$\dd_h= \dd _0  |h|/t,\ \ \ \gg_h(\dd,\x)= \ff{c|h|^2}\dd \Big(\ff 1 {t^3} +\ff{\log\tt W(\x)}{t^2}\Big).$$ Therefore, the desired Harnack inequality follows from Proposition \ref{P4.0}.
\end{proof}

To derive the log-Harnack inequality, we need the following slightly stronger condition than the second one in {\bf (A)}: there exists an increasing function $U$ on $[0,\infty)$ such that
\beq\label{WFY} |Z(\x)-Z(\y)|^2\le |\x-\y|^2 \big\{U(|\x-\y|) + \ll W(\y)\big\},\ \ \x,\y\in \R^{m+d}.\end{equation}

\beg{thm}\label{T4.1} Assume {\bf (A)} such that $(\ref{WFY})$ holds. Then there exists a constant $c>0$ such that
$$P_t\log f(\x) - \log P_t f(\y)
 \le c|\x-\y|^2\bigg\{\ff 1 {(1\land t)^3} +  \ff {U((1\lor t^{-1})|\x-\y|) +   W(\y)}{t\land 1}\bigg\}$$ holds for any $t>0,$ positive function $f\in \B_b(\R^{m+d})$, and $\x,\y \in\R^{m+d}.$
\end{thm}

\beg{proof} Again as in the proof of Corollary \ref{C2}, it suffices to  prove  for $t\in (0,1].$ Let $\x=(x,y)$ and $\y=(\tt x,\tt y)$. We will make use of the coupling constructed in the proof of Theorem \ref{T1.1} for $\vv=1, h= (x-\tt x, y-\tt y) $ and $(u,v)$ being in (\ref{uv1}). We have $(X_t,Y_t)=(X_t^1,Y_t^1)$, and $(X_s^1,Y_s^1)_{s\in [0,t]}$ is generated by $L$ under the probability $\Q_1= R_t^1\P$. So, by the Young inequality (see \cite[Lemma 2.4]{ATW09}), we have
 \beg{equation*}\beg{split} P_t \log f(\tt x,\tt y) &=\E \big(R_t^1\log f(X_t^1, Y_t^1)\big) = \E\big(R_t^1\log f(X_t,Y_t)\big)\\
 &\le \E(R_t^1\log R_t^1)+ \log\E f(X_t,Y_t)= \log P_t f(x,y)+ \E(R_t^1\log R_t^1).\end{split}\end{equation*} Combining this with (\ref{NNm}) we arrive at
\beq\label{BC} P_t\log f(\tt x,\tt y) - \log P_t f(x,y)\le \ff 1 2 \E_{\Q_1}\int_0^t |\si^{-1}\xi_s^1|^2\d s.\end{equation}
Taking $z$ such that $|z|\le \|A^{-1}\|\cdot |h_1|$, we obtain from   (\ref{2.2'}), (\ref{WFY}), (\ref{uv1}) and (\ref{uv2}) that
\beg{equation*}\beg{split} |\si^{-1}\xi_s^1|^2 & \le  \Big\{|\LL(h,z,s)|^2 +|\Theta(h,z,s)|^2\big(U(|\Theta(h,z,s)|)+\ll W(X_s^1,Y_s^1)\big)\Big\}\\
&\le c|h|^2\Big\{\ff 1 {t^4} +  \ff {  U(|h|/t) +   W(X_s^1,Y_s^2)}{t^2}\Big\}.\end{split}\end{equation*}
Combining this with (\ref{BC}) and noting that $LW\le CW$ implies $\E_{\Q_1} W(X_s^1, Y_s^1)\le \e^{Cs}W(\tt x,\tt y)$ for $s\in [0,t]$, we complete the proof.
\end{proof}

We conclude this part, we come back to Example 2.1 for the kinetic Fokker-Planck equation.

\paragraph{Example 4.1} In Example 2.1 let e.g. $V(x)= (1+|x|^2)^l$. Then {\bf (A)} and (\ref{WFY}) holds for $W(x,y)= \exp[2V(x)+|y|^2]$ and
$U(r)= c r^{2[(2l-1)\lor 1]}$ for some constant $c>0.$ Therefore, Theorem \ref{T4.1} applies.

Next, for the gradient-entropy inequality (\ref{GGA}) and (\ref{HHH}), let us consider for simplicity that $m=d=1$ and $V(x)=x^3$:
\beq\label{kfp2} \beg{cases} \d X_t= Y_t\d t,\\
\d Y_t= \d B_t -(X_t)^3\d t-Y_t\d t.\end{cases}\end{equation}
In this case we have $Z(x,y)=-x^3-y$,  so that
$$|Z(x,y)-Z(\tt x,\tt y)|^2\le c(|x-\tt x|^2+|y-\tt y|^2)(1+x^4+\tt x^4).$$ Next, let $W(x,y)=1+ \ff 1 2 x^4+ y^2$. We have
$$LW(x,y)= 2y x^3 + 1 -2 x^3 y -2 y^2 = 1- 2 y^2 \le W(x,y).$$ Thus, (\ref{A1}) holds for $U=0$. Moreover,
following the line of  in  \cite{Wu01,BCG08,DFG09}, consider $w(x,y)=a\big(\ff 1 2 x^4+y^2)+bxy$ for some well chosen constant $a,b$ and putting ${\tt W}(x,y)=\exp(w-\inf w)$,  we have
$$-\ff{L{\tt W}}{\tt W}\ge \aa W-K$$ for some constants $\aa,K>0$.
Indeed,
\beg{equation*}\beg{split}
\ff{L \tilde W}{\tt W}(x,y) &= L\log \tt W(x,y)- \ff 1 2 |\pp_y \log \tt W|^2(x,y)\\
&= a+2a^2y^2-2ax^3y-bx^4-2ay^2-bxy+2ax^3y+by\\
&\le a+(2a^2-2a+b(1+\vv/2))y^2-bx^4+bx^2/(2\vv)\\
&\le K-\aa (1+y^2+x^4)\end{split}\end{equation*}
holds for some constants $\aa,K>0$ by taking $a,b,\vv>0$ such that   $2a^2-2a+b(1+\varepsilon/2)<0$.
Therefore, (\ref{A2}) holds for some $\ll, K>0$ so that    (\ref{GGA}) and (\ref{HHH}) hold.

\subsection{Harnack inequality under assumption (H)}

As shown in \cite{ATW09}, the  derivative estimate  (\ref{1.4}) will enable us to prove an Harnack inequality with a power in the sense of \cite{W97}.  More precisely, we have the following result.

\beg{thm} \label{THarnack}   Let  $(\ref{NZ})$  hold and let $\Phi_t$ be in $(\ref{Phi})$. Then for any $t>0,\aa>1$ and positive function $f\in \B_b(\R^{m+d}),$
\beg{equation} \label{Har}(P_t f)^\aa(x,y)\le (P_t f^\aa)(\tt x,\tt y)\exp\Big[\ff{\aa}{\aa-1}\Phi_t(|x-\tt x|,|y-\tt y|)\Big],\ \ (x,y), (\tt x,\tt y)\in \R^{m+d}\end{equation}
holds.    Consequently,
 \beq\label{LHar} P_t \log f(x,y)\le \log P_t f(\tt x,\tt y) + \Phi_t(|x-\tt x|,|y-\tt y|),\ \  (x,y), (\tt x,\tt y)\in \R^{m+d}.\end{equation}
\end{thm}
\beg{proof} It is easy to see that (\ref{Har}) follows from (\ref{1.4}) and Proposition \ref{P4.0}. Next, according to \cite[Proposition 2.2]{W10},  (\ref{LHar}) follows from $(\ref{Har})$ since $\R^{m+d}$ is a length space under the metric
$$\rr((x,y),(\tt x,\tt y)):= \ss{\Phi_t(|x-\tt x|, |y-\tt y|)}.$$ So,   (\ref{Har}) implies (\ref{LHar}).
\end{proof}

The next result extends Theorem \ref{THarnack} to unbounded $\nn Z$.

\beg{thm} \label{THarnack2}   Assume {\bf (H)}.  Then there exists a constant $c>0$ such that for any $t >0$ and positive $f\in\B_b(\R^{m+d})$,
\beq\label{LHar'} \beg{split} & P_t \log  f(\y)- \log P_t f(\x)\\
  & \le |\x-\y|^2  \Big\{\ff c {(1\land t)^3}
+  \ff{  c}{(1\land t)^{2l}}\big(1+|\x|+|\y|\big)^{4l}\Big\}\end{split}\end{equation} holds for  $ \x,\y\in \R^{m+d}.$
  If {\bf (H)} holds for some $l<\ff 1 2$, then there exists a constant $c>0$ such that
\beg{equation} \label{Har'} \beg{split}&(P_t f)^\aa(\x)\le (P_t f^\aa)(\y) \\
&\times \exp\Big[\ff{\aa c|\x-\y|^2}{(\aa-1)(1\land t)^4}
\Big\{(|\x|\lor |\y|)^{4l} + \big((\aa-1)(1\land t)^2\big)^{4(l-1)/(1-2l)}\Big\}\Big]\end{split}\end{equation}
holds for all $t>0,\aa>1,   \x,\y\in \R^{m+d}$ and positive $f\in \B_b(\R^{m+d}).$
\end{thm}

\beg{proof}  (\ref{LHar'}) follows from Theorem \ref{T4.1} since in this case {\bf (A)} and (\ref{WFY}) hold for $W(\x)= (1+|\x|^2)^{2l} $ and $U(r)=c r^{2l}$ for some $\ll,c>0;$ while  (\ref{Har'}) follows from Corollary \ref{C3}(2) and Proposition \ref{P4.0}.
 \end{proof}

According to \cite[Proposition 2.4]{W10}, we have   the following   consequence of Theorems \ref{THarnack} and \ref{THarnack2}.

\beg{cor} \label{C1.4} Let $p_t$ be the transition density of $P_t$ w.r.t. some $\si$-finite  measure $\mu$ equivalent to the Lebesgue measure on $\R^{m+d}$.   Let $\Phi_t$ be in Theorem $\ref{THarnack}$.
 \beg{enumerate} \item[$(1)$]   $(\ref{NZ})$  implies
\beg{equation*}\beg{split} &\int_{\R^{m+d}} \bigg(\ff{p_t((x,y),\z)}{p_t((\tt x,\tt y),\z)} \bigg)^{1/(\aa-1)} p_t((x,y),\z)
\mu(\d \z)\le\exp\bigg[\ff{\aa }{(\aa-1)^2}\Phi_t(|\tt x-x|, |\tt y-y|)\bigg],\\
&\int_{\R^{m+d}} p_t((x,y),\z)\log\ff{p_t((x,y),\z)}{p_t((\tt x,\tt y),\z)} \mu(\d\z)\le
\Phi_t(|\tt x-x|, |\tt y-y|).\end{split}\end{equation*}    for any $t>0$ and $(x,y),(\tt x,\tt y)\in \R^{m+d}.$
\item[$(2)$] If {\bf (H)} holds for some $l\in (0,\ff 1 2)$, then there exists a constant $c>0$ such that
\beg{equation*}\beg{split} &\int_{\R^{m+d}} \bigg(\ff{p_t(\x,\z)}{p_t(\y,\z)} \bigg)^{1/(\aa-1)} p_t(\x,\z)\mu(\d \z)\\
&\le\exp\bigg[\ff{\aa c|\x-\y|^2}{(\aa-1)^2(1\land t)^4}
\Big\{(|\x|\lor |\y)^{4l} + \big((\aa-1)(1\land t)^2\big)^{4(l-1)/(1-2l)}\Big\}\bigg]\end{split}\end{equation*}
 holds  for all $t>0$ and $\x,\y\in \R^{m+d}.$
\item[$(3)$] If {\bf (H)} holds then there exists a constant $c>0$ such that
\beg{equation*}\beg{split} &\int_{\R^{m+d}} p_t(\x,\z)\log\ff{p_t(\x,\z)}{p_t(\y,\z)} \mu(\d\z)\\
&\le
|\x-\y|^2  \Big\{\ff c {(1\land t)^3}
+  \ff{c}{(1\land t)^{2l}}\big(1+|\x|+|\y|\big)^{4l}\Big\}\end{split}\end{equation*}  holds  for all $t>0$ and $\x,\y\in \R^{m+d}.$ \end{enumerate} \end{cor}

Next, for two probability measures $\mu$ and $\nu$, let
$\C(\nu,\mu)$ be the class of their couplings, i.e. $\pi\in \C(\nu,\mu)$ if $\pi$ is a probability meadsure on $\R^{m+d}\times \R^{m+d}$ such that $\pi(\R^{m+d}\times \cdot)=\mu(\cdot)$ and $\pi(\cdot\times \R^{m+d})=\nu(\cdot)$.
Then according to the proof of \cite[Corollary 1.2(3)]{RW}, Theorems \ref{THarnack} and \ref{THarnack2} also imply the following entropy-cost inequalities. Recall that for any non-negative symmetric  measurable function ${\bf c}$ on $\R^{m+d}\times \R^{m+d}$, and for any two probability measures $\mu,\nu$ on $\R^{m+d}$, we call
$$W_{\bf c}(\nu,\mu):= \inf_{\pi\in\scr C(\nu,\mu)}\int_{\R^{m+d}\times\R^{m+d}}{\bf c}(\x,\y)\,\d\pi(\d \x,\d \y)$$ the transportation-cost between these two distributions induced by the cost function ${\bf c}$, where $\scr C(\nu,\mu)$ is the set of all couplings of $\nu$ and $\mu$.

\beg{cor} \label{C4b} Let $P_t$ have an invariant probability measure $\mu$, and let $P^*$ be the adjoint operator of $P$ in $L^2(\mu).$ \beg{enumerate} \item[$(1)$] If   $(\ref{NZ})$ holds then
\beq\label{TC} \mu(P_t^*f\log P_t^*f)\le    W_{{\bf c}_t}(f\mu,\mu),\ \ t>0, f\ge 0,\mu(f)=1,\end{equation} where
${\bf c}_t(x,y;\tt x,\tt y)= \Phi_t(|\tt x-x|, |\tt y-y|).$
\item[$(2)$]  If {\bf (H)} holds, then there exists $c>0$ such that
$(\ref{TC})$ holds for
$${\bf c}_t(\x, \y)=  |\x-\y|^2  \Big\{\ff c {(1\land t)^3}
+  \ff{ c}{(1\land t)^{2l}}\big(1+|\x|+|\y|\big)^{4l}\Big\}.$$\end{enumerate}\end{cor}

\paragraph{Remark 4.1}
\begin{enumerate}
\item[(I)] Recall that the Pinsker inequality says that for any two probability measures $\mu,\nu$ on a measurable space, the total variation norm of $u-v$ is dominated by the square root of   twice relative entropy of $\nu$ w.r.t. $\mu$. Combining this inequality with  (1) of Corollary \ref{C1.4}, assuming thus $\|\nabla Z\|_\infty<\infty$, we get
$$\|P_t((x,y),\cdot)-P_t((\tt x,\tt y),\cdot)\|_{TV}\le \sqrt{2\,\Phi_t(|\tt x-x|, |\tt y-y|)},$$
which may be useful as an alternative to small set evaluation in Meyn-Tweedie's approach for convergence to equilibrium for the kinetic Fokker-Planck equation.
\item[(II)] Using Villani's result \cite[Th.39]{Vil09} in the kinetic Fokker Planck case which asserts that if $|\nabla^2V|$ is bounded and $\mu$ as a product measure satisfies a logarithmic Sobolev inequality, then there is an exponential convergence towards equilibrium in entropy, so that
$$\mu(P_s^*f\log P_s^*f)\le Ce^{-Ks}\mu(P_1^*f\log P_1^*f),\ \ f\ge 0, \mu(f)=1, s\ge 1$$
holds for some constant $C>0$. Combining this with Talagrand inequality implied by the logarithmic Sobolev inequality (see \cite{OV00}) and using   Corollary \ref{C4b}(1),  we get
$$W_2^2(P_sf\mu,\mu)\le C'e^{-Ks}W_2^2(f\mu,\mu),\ \ f\ge 0, \mu(f)=1, s\ge 1$$ for some constant $C'>0$,
where $W_2^2= W_{\bf c}$ for ${\bf c}(\x;\y):= |\x-\y|^2.$ This
  generalizes the exponential convergence in Wasserstein distance derived in \cite{BGM10} for the non interacting case.
\end{enumerate}

\paragraph{Acknowledgement} The second named author would like to thank Professor Xicheng Zhang for introducing his very interesting   paper
\cite{Zhang}. Both authors would like to thank the referee for helpful comments.

\beg{thebibliography}{99}

\bibitem{AT99}  M. Arnaudon, A. Thalmaier, \emph{Bismut type differentiation of semigroups,} Probability Theory and Mathematical Statistics 23--32, VSP/TEV, Utvecht and Viluius, 1999.

\bibitem{ATW06} M. Arnaudon, A. Thalmaier, F.-Y. Wang,
  \emph{Harnack inequality and heat kernel estimates
  on manifolds with curvature unbounded below,} Bull. Sci. Math. 130(2006), 223--233.

\bibitem{ATW09} M. Arnaudon, A. Thalmaier, F.-Y. Wang,
  \emph{Gradient estimates and Harnack inequalities on non-compact Riemannian manifolds,}
   Stoch. Proc. Appl. 119(2009), 3653--3670.

\bibitem{BCG08} D. Bakry, P. Cattiaux, A. Guillin, \emph{Rate of convergence for ergodic continuous Markov processes : Lyapunov versus Poincare,}
J. Func. Anal. 254 (2008), 727--759.

\bibitem{Bismut} J. M. Bismut, \emph{Large Deviations and the
Malliavin Calculus,} Boston: Birkh\"auser, MA, 1984.

\bibitem{BGM10} F. Bolley, A. Guillin, F. Malrieu, \emph{Trend to equlibrium and particle approximation for a weakly selfconsistent Vlasov-Fokker-Planck equation}, M2AN 44(5) (2010), 867--884.

\bibitem{DRW09}  G. Da Prato, M. R\"ockner, F.-Y. Wang, \emph{Singular stochastic equations on Hilbert
spaces: Harnack inequalities for their transition semigroups,} J. Funct. Anal. 257 (2009), 992--017.

\bibitem{DFG09} R. Douc, G. Fort, A. Guillin, \emph{Subgeometric rates of convergence of f-ergodic strong Markov processes},
Stoch. Proc. Appl. 119 (2009) 897--923.

\bibitem{EL} K.D. Elworthy,  Xue-Mei Li, \emph{Formulae for the
derivatives of heat semigroups,} J. Funct. Anal. 125(1994),
252--286.

\bibitem{ES} A. Es-Sarhir, M.-K. v. Renesse, M. Scheutzow,
\emph{Harnack inequality for functional SDEs with bounded memory,}
Electron. Commun. Probab. 14 (2009), 560--565.

\bibitem{K}  H. Kawabi,  \emph{The parabolic Harnack inequality for
the time dependent Ginzburg-Landau type SPDE and its application,}
Pot. Anal.  22(2005), 61--84.

\bibitem{LW} W. Liu, F.-Y. Wang, \emph{Harnack inequality and strong Feller
  property for stochastic fast diffusion equations,} J. Math. Anal. Appl.
  342(2008), 651--662.

\bibitem{OV00} F. Otto, C. Villani, \emph{Generalization of an inequality by Talagrand and links with the logarithmic
Sobolev inequality}, J. Funct. Anal., 173 (2000), 361--400.

\bibitem{Ouyang} S.-X. Ouyang, \emph{Harnack inequalities and applications for multivalued stochastic evolution equations,}   Infin. Dimens. Anal. Quant. Probab.  Relat. Topics 14(2011), 261--278.

\bibitem{ORW} S.-X. Ouyang, M. R\"ockner, F.-Y. Wang,
 \emph{Harnack inequalities and applications for Ornstein-Uhlenbeck semigroups with
 jump,}   Pot. Anal.   36(2012), 301--315.

\bibitem{RW} M. R\"ockner, F.-Y. Wang, \emph{Log-Harnack  inequality for stochastic differential equations in Hilbert spaces and its consequences, } Infin. Dimens. Anal. Quant. Probab.  Relat. Topics 13(2010), 27--37.

 \bibitem{Vil09} C. Villani, \emph{Hypocoercivity}.  Mem. Amer. Math. Soc.  202  (2009),  no. 950.

\bibitem{W97} F.-Y. Wang,  \emph{Logarithmic Sobolev
inequalities on noncompact Riemannian manifolds,} Probability
Theory Relat. Fields 109(1997), 417--424.

\bibitem{W07} F.-Y. Wang, \emph{Harnack inequality and applications
for stochastic generalized porous media equations,}  Ann. Probab. 35(2007),
1333--1350.

\bibitem{W10} F.-Y. Wang,  \emph{Harnack inequalities on manifolds with boundary and applications,}  J. Math. Pures Appl. 94(2010), 304--321.

\bibitem{W10b} F.-Y. Wang, \emph{Coupling and its applications}, Preprint, accessible on arXiv:1012.5687.

\bibitem{WX10} F.-Y. Wang, L. Xu, \emph{Derivative formula and applications for hyperdissipative stochastic Navier-Stokes/Burgers equations,} to appear in Infin. Dimens. Anal. Quant. Probab.  Relat. Topics,  accessible on  arXiv:1009.1464.

\bibitem{WWX10}  F.-Y. Wang, J.-L. Wu, L. Xu, \emph{Log-Harnack inequality for stochastic Burgers equations and applications,}   J. Math. Anal. Appl. 384(2011), 151--159.

\bibitem{WY10} F.-Y. Wang, C. Yuan, \emph{Harnack inequalities for functional SDEs with multiplicative noise and applications}, Stoch. Proc. Appl. 121(2011), 2692--1710.

\bibitem{Wu01} L. Wu, \emph{Large and moderate deviations and exponential convergence for stochastic damping Hamiltonian systems}, Stoch. Proc. Appl., 91 (2001), 205--238.

\bibitem{Z} T.-S. Zhang, \emph{White noise driven SPDEs with reflection: strong Feller properties and Harnack inequalities,} Potential Anal.
 33 (2010),137--151.

\bibitem{Zhang} X. Cheng, \emph{Stochastic flows and Bismut formulas for stochastic Hamiltonian systems,} Stoch. Proc. Appl. 120(2010), 1929--1949.

 \end{thebibliography}
\end{document}